\documentclass{amsart}
\usepackage{amsmath}
\usepackage{amssymb}
\usepackage{graphicx}
\usepackage{url}

\begin{document}

\title[Calling a spade a spade]{Calling a spade a spade:\\ Mathematics in the new pattern of division of labour}
\author{Alexandre V. Borovik}

\email{alexandre$\gg$at$\ll$borovik.net}
\thanks{The last pre-publication version, 11 December 2014. \copyright 2014 Alexandre Borovik}

\newtheorem{ex}{Problem}
\newcommand{\bq}{\begin{quote}}
\newcommand{\eq}{\end{quote}}
\newcommand{\bi}{\begin{itemize}}
\newcommand{\ei}{\end{itemize}}
\newcommand{\fff}[1]{\medskip\noindent {{\bf #1.}}}

\maketitle

\small
\begin{flushright}
\emph{ The man who could call a spade a spade\\ should be compelled to use one.\\ It is the only thing he is fit for.}\\

Oscar Wilde\\

\end{flushright}\normalsize

\bigskip\noindent

\section{Introduction}

The growing disconnection of the majority of the population from mathematics is increasingly difficult to ignore.

This paper focuses on the socio-economic roots of this cultural and social phenomenon which are not usually mentioned in public debates. I concentrate on mathematics education, as an important and well documented area of interaction of mathematics with the rest of human culture.

New patterns of division of labour have dramatically changed the nature and role of mathematical skills needed for the labour force and correspondingly changed the place of mathematics in popular culture and in mainstream education. The forces that drive these changes come from the tension between the ever deepening specialisation of labour and ever increasing length of specialised learning required for jobs at the increasingly sharp cutting edge of technology.

Unfortunately these deeper socio-economic origins of the current systemic crisis of mathematics education are not clearly spelt out, neither  in  cultural studies nor, even more worryingly, in the education policy discourse;  at the best, they are only euphemistically hinted at.

This paper is an attempt to describe the socio-economic landscape of mathematics education without resorting to euphemisms. This task imposes on the author certain restrictions: he cannot take sides in the debate and therefore has to refrain from giving any practical recommendations. Also it makes necessary a very clear disclaimer:
\bq
\small
\emph{The author writes in his personal capacity. The views expressed do not necessarily represent the position of his employer or any other person, organisation, or institution.}
\normalsize
\eq

\newpage
\section{The new division of labour}

\small
\begin{flushright}
\emph{It's the economy, stupid.}\\

James Carville\footnote{\emph{It's the economy, stupid.} According to Wikipedia, this phrase, frequently attributed to Bill Clinton, was made popular by James Carville, the strategist of Clinton's successful 1992 presidential campaign against George H. W. Bush.}\\
\end{flushright}
\normalsize

\subsection{A word of wisdom from Adam Smith}

Discussion of mathematics education takes place in a  socioeconomic landscape which has never before existed in the history of humanity.

This, largely unacknowledged, change, can be best explained by invoking  Adam Smith's famous words displayed on the British \pounds 20 banknote, Figure~\ref{fig:20pounds}:

\begin{figure}[h]
\begin{center}
\includegraphics[width=4.5in]{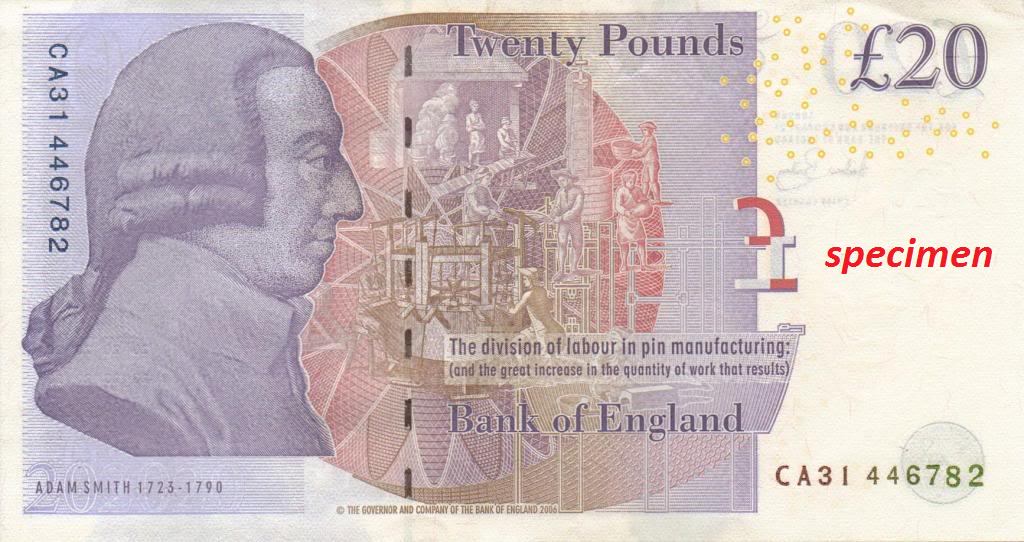}
\end{center}
\caption{}
\label{fig:20pounds}
\end{figure}

The words on the banknote:
\bq
\small
The division of labour in pin manufacturing  (and the great increase in the quantity of work that results)
\normalsize
 \eq
are, of course, a quote from Adam Smith's \emph{The Wealth of Nations}. They are found on the very first page of
Chapter I of Book I with the now famous title \emph{Of The Division of Labour}:

\begin{quote}
\small
One man draws out the wire; another straights it; a third cuts
it; a fourth points it; a fifth grinds it at the top for receiving
the head; to make the head requires two or three distinct
operations; to put it on is a peculiar business; to whiten the
pins is another; it is even a trade by itself to put them into the
paper; and the important business of making a pin is, in this
manner, divided into about eighteen distinct operations.
\normalsize
\end{quote}

And Adam Smith comes to the conclusion:

\begin{quote}
\small
\dots\ they certainly could not each of them have made twenty, perhaps
not one pin in a day; that is, certainly, not the two hundred and
fortieth, perhaps not the four thousand eight hundredth part of what
they are at present capable of performing, \dots
\normalsize
\end{quote}

By the start of the 21st century, the ever deepening division of labour has reached a unique point in the history of humankind when 99\% of people have not even the vaguest idea about the workings of 99\% of technology in their immediate surrounding---and this applies even more strongly to technological uses
of mathematics, which are mostly invisible.

Every time you listen to an iPod or download a compressed graphic file from the Internet, extremely sophisticated mathematical algorithms come into play. A smartphone user never notices this because these algorithms are encoded deep inside the executable files of smartphone apps. Nowadays mathematics (including many traditional areas of abstract pure mathematics, such as number theory, abstract algebra, combinatorics, and spectral analysis, to name a few) is used in our everyday life thousands, maybe millions, of times more intensively than  50 or even 10 years ago. Mathematical results and concepts involved in practical applications are much deeper and more abstract and  difficult than ever before. One of the paradoxes of modern times is that this makes mathematics invisible because it is carefully hidden behind a user friendly smartphone interface.

There are more mobile phones  in the world now than toothbrushes. But the mathematics built into a mobile phone or an MP3 player is beyond the understanding of most graduates from mathematics departments of British universities. However, practical necessity forces us to teach a rudimentary  MP3/MP4 technology, in cookbook form, to electronic engineering students; its  mathematical content is diluted or even completely erased.

\subsection{A few more examples}

New patterns of division of labour manifest themselves at every level of the economy.

\subsubsection{A consumer} 25 years ago in the West, the benchmark of arithmetic competence at the consumer level was the ability to balance a chequebook. Nowadays, bank customers can instantly get full information about the state of their accounts from an app on a mobile phone.

\subsubsection{A worker in the service sector} How much arithmetic should a worker at a supermarket checkout know? And they are being replaced by fully automated self-checkout machines.

\subsubsection{A worker in an old industry} Even in the pre-computer era, say, in the 19th and the first half of 20th centuries consumers were increasingly ignorant of the full extent of technological sophistication used in the production of everyday goods. In relation to mathematics that meant that buyers of ready-to-wear clothing, for example, were likely to be unaware of craft-specific shortcuts and tricks of geometry and arithmetic used by a master cutter when he made a template for a piece of clothing. In the clothing industry nowadays, cutters are replaced by laser cutting machines. But a shirt remains essentially the same shirt as two centuries ago; given modern materials, a cutter and a seamstress of yesteryear would still be able to produce a shirt meeting modern standards (and millions of seamstresses are still toiling in the  sweatshops of the Third World). What a 19th or 20th century cutter would definitely not be able to do is to develop mathematical algorithms which, after being converted into computer code, control a laser cutting machine. Design and optimisation of these algorithms require a much higher level of mathematical skills and are mostly beyond the grasp of the majority of our mathematics graduates.

\subsubsection{A worker in a new industry} Do you need any mathematical skills at all for snapping mobile phones together on an assembly line? But production of microchips is highly automated and involves a very small number of much better trained and educated workers. Research and development in the area of microelectronics (and photonics) is of course an even more extreme case of concentration of expertise and skills.

\subsubsection{International division of labour} It is easy to imagine a country  where not a single person has a working knowledge of semiconductor technology and production of microchips. What for? Microchips are something sitting deep inside electronic goods imported from China---and who cares what is inside? Modern electronic goods usually have sealed shells, they are not supposed to be opened. Similarly, one can easily imagine a fully functioning country where no-one has mastered, say, long division or factorisation of polynomials.

\subsection{Social division of labour}

In the emerging division of {intellectual labour}, {mathematics is a 21st century equivalent of sharpening a pin.}

The only difference is that a pin-sharpener of Adam Smith's times could be trained on the job in a day. Development of a mathematically competent worker for high tech industries  requires at least 15 years of education from ages  $5$ to $20$.

\bq
\textbf{It is this tension between the ever-increasing degree of specialisation and the ever-increasing length of specialised education that lies at the heart of the matter.}
\eq

At this point we need to take a closer look at \emph{social division} of labour. Braverman \cite{Braverman74} emphasises the distinction between the \emph{social} division of labour between different occupational strata of society and the \emph{detailed} division of labour between individual workers in the workplace.

\bq\small
The division of labor in society is characteristic of all known societies; the division of labor in the workshop is the special product of capitalist society. The social division of labor divides society among occupations, each adequate to a branch of production; the detailed division of labor destroys occupations considered in this sense, and renders the worker inadequate to carry through any complete production process. In capitalism, the social division of labor is enforced chaotically and anarchically by the market, while the workshop division of labor is imposed by planning and control.  \cite[pp.~50--51]{Braverman74}
\normalsize\eq

It is the new workplace, or ``detailed'', division of labour that makes mathematics redundant in increasingly wide areas of professional occupation. Meanwhile the length-of-education constraints in reproduction of a mathematically skilled workforce lead to mathematics being single out not only in workplace division of labour, but also in social division. And, exploiting the above quote from Braverman, it is the  ``chaotic and anarchic'' nature of social division  that leads to political infighting around mathematics education and paralyses education policy making.

The rest of my paper expands on these  theses. One point that I do not mention is the division of labour \emph{within} mathematics; this is an exciting topic, but it requires a much more specialised discussion.

\section{Politics and economics}

The issue of new patterns of division of labour has begun to emerge in political discourse. I give here some examples.

The book by Frank Levy and Richard Murnane \emph{The New Division of Labor} \cite{Levy04}, published in 2004 and based on  material from the USA, focuses on economic issues viewed from a business-centred viewpoint. Here is a characteristic quote:

\bq\small
In economic terms, improved education is required to restore the labor market to balance. [\dots ] the falling wages of lower skilled jobs reflect the fact that demand was not keeping up with supply. If our predictions are right, this trend will continue as blue-collar and clerical jobs continue to disappear.

Better education is an imperfect tool for this problem. The job market is changing fast and improving education is a slow and difficult process. \cite[p.~155]{Levy04}.

\normalsize\eq

Elizabeth Truss, a Conservative Member of Parliament and Secretary for the Environment (who until recently was an Undersecretary of State in the Department for Education), not long ago published a report \cite{Truss11} where she addressed the issue of the ``hourglass economy'' in the context of education policy.

\bq\small
The evidence suggests
increased polarisation between high
skilled and unskilled jobs, with skilled
trades and clerical roles diminishing.
Long standing industries are becoming
automated, while newly emerging
industries demand high skills. Formal
and general qualifications are the main
route into these jobs. At the top level
MBAs and international experience is
the new benchmark. Despite popular
perception, the middle is gradually
disappearing to create an `hourglass
economy'. \cite[p.~1]{Truss11}
\normalsize\eq

In the next section, we shall return to the ``hourglass economy'' and the ``hourglass'' shape of the demand for mathematics education to different levels of students' attainment. Meanwhile, I refer the reader to the views of numerous economists concerning ``job polarisation'' (Autor \cite{Autor10}, Goos et al. \cite{Goos09}), ``shrinking middle'' (Abel and Deitz \cite{abel-deitz}), ``intermediate occupations'' and  ``hourglass economy'' (Anderson \cite{Anderson09}).
The same sentiments about the ``disappearing middle'' are repeated in more recent books under catchy titles such as Tyler Cowen's \emph{The Average is Over} \cite{Cowen13}; they are becoming part of the \emph{Zeitgeist}. Although their book is optimistic, Brynjolfsson and  McAfee \cite{Brynjolfsson14}  emphasise the way in which the application of the know-how in the upper half of the hourglass causes the hollowing out of  the ``neck''.

It is instructive to compare opinions on job polarisation and its impact on education coming from opposite ends of the political spectrum.

Judging by his recent book \cite[Chapter 14]{Greenspan13}, Alan Greenspan focuses on the top part of the hourglass:

\bq\small
[W]e may not have the capability to educate and train students up to the level increasingly required by technology to staff our ever more complex capital stock. The median attainment of our students just prior to World War II was a high school diploma. That level of education at the time, with its emphasis on practical shop skills, matched the qualifications, by 1950s standards, for a reasonably skilled job in a steel mill or auto-assembly plant. [\dots ] These were the middle income jobs of that era.

\textbf{But the skill level required of today's highly computerized technologies cannot be adequately staffed with today's median skills.} [The emphasis is mine---AB.]

\normalsize
\eq

A voice from the left (Elliot \cite{Elliott11}), on the contrary, suggests that education has been intentionally dumbed down:

\bq
\small
We need, I should say, to look for an analysis in the direction of global developments in the capitalist labour process---especially the fragmentation of tasks, the externalization of knowledge (out of human heads, into computer systems, administrative systems and the like)---and the consequent declining need, among most of the population, regarded as employees or workers, for the kinds of skills (language skills, mathematical skills, problem-solving skills etc.) which used to be common in the working class, let alone the middle classes. This analysis applies to universities and their students. Dumbing-down is a rational---from the capitalist point of view---reaction to these labour-process developments. No executive committee of the ruling class spends cash on a production process (the production of students-with-a-diploma) that, from its point of view, is providing luxury quality. It will continuously cut that quality to the necessary bone. It is doing so. This, to repeat the point, is a global tendency rooted in the reality of capitalist production relations.
\normalsize
\eq

But Greenspan \cite[Chapter 14]{Greenspan13} appears to take a more relaxed   view on changes in economic demand for education:
\bq
\small
{While there is an upside limit to the average intellectual capabilities of population, there is no upper limit to the complexity of technology.
\normalsize}
\eq
\bq
\small
{With [\dots ] an apparently inbred upper limit to human IQ, are we destined to have an ever smaller share of our workforce staff our ever more sophisticated high-tech equipment and software?
\normalsize}
\eq

 Many  may disagree with this  claim---but it may nevertheless influence political and business decisions.

\section{Implications for mathematics education}

We have to realise that it is no longer an issue whether the role of  mathematics in society is changing: the change is being ruthlessly  forced on us by Adam Smith's `invisible hand'.

\begin{figure}[h]
\begin{center}
\includegraphics[width=4.5in]{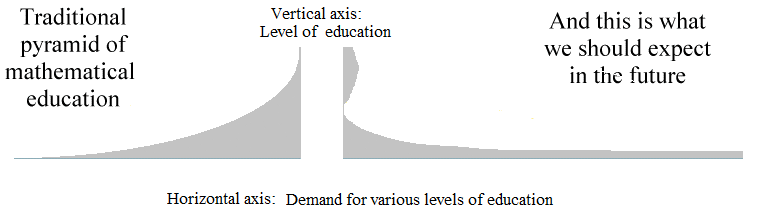}
\end{center}
\caption{Pyramids of economic demand for mathematics education (qualitative schemes, not to scale, but higher levels of education correspond to higher levels in the pyramids).}
\label{fig:piramids}
\end{figure}

In particular, changing economic imperatives lead to the collapse of the traditional pyramid of mathematics education. Let us look at the diagram in Figure~\ref{fig:piramids}.

The diagram is not made to any scale and should be treated qualitatively, not quantitatively.
The left hand side of the pyramid suggests how the distribution of mathematical attainment looked in the mid 20th century, with pupils / students / graduate students at every level of education being selected from a larger pool of students at the previous level. In the not so distant past, every stage in mathematics education matched the economic demand for workers with a corresponding level of skills. From  students' point of view, every  year invested in mathematics education was bringing them a potential (and immediately cashable) financial return.

The traditional pyramid of mathematics education was stable because every level had its own economic justification and employment opportunities. I have included as the Appendix  the \emph{Post Office Entrance Examination} from 1897 which is being circulated among British mathematicians as a kind of subversive leaflet. A century ago, good skills in practical arithmetic opened up employment opportunities for those in the reasonably wide band of the diagram on the left, the one which has now become the bottleneck of the `hourglass' on the right. Nowadays this level of skills is economically redundant; its only purpose is to serve as an indication of, and as a basis for, a person's progress to higher, more economically viable, levels of mathematics education.

The right hand side of the pyramid suggests what we should expect in the future: an hourglass shape, with intermediate levels eroded. Certain levels of mathematics education are not supported by immediate economic demand and serve only as an intermediate or preparatory step for further study. From an individual's point of view, the economic return on investment in mathematical competence is both delayed and less certain. Once this is realised, it seems likely to weaken the economic motivation for further study.

Many practitioners of mathematics education \cite{Edwards14} and sociologists \cite{Gainsburg05} are coming to the same conclusion:
\bq
\small
Studies
of the actual demands of everyday adult practices reveal that most occupations involve only a low
level of mathematical content and expose the disparate natures of everyday and school mathematics. \cite[p.~1]{Gainsburg05}
\Large
\eq
\bq
\small
 [\dots ] most jobs that currently require advanced technical degrees are using that requirement simply as a filter. \cite[p.~21]{Edwards14}
 \Large
\eq

The cumulative nature of learning mathematics makes a ``top-heavy'' model of education \textbf{unsustainable}: what will be the motivation for students to struggle through the  neck of the hourglass? Whether they realise it or not (most likely not) children and their families  subconsciously apply a discounted cash flow analysis to the required intellectual effort and investment of time as compared to the subsequent reward.

Education (or at least state-run education) is a sector of the economy where real consumer choice does not exist. Of course, there are a couple of choice points at which students and their families can decide what to study---but not how. There is no real choice of schools and teachers. From the economics point of view, the state education system in England is the same as the state education in the former communist block (and this phrase is not intended as criticism of either of them).
But it is the \textsc{Aeroflot} business model of yesteryear:

\bq
\small
\textsc{Aeroflot} flight attendant:
\bq
\small
``\emph{Would you like a dinner}?''
\eq
Passenger:
\bq
\small
``\emph{And what's the choice}?''
\eq
Flight attendant:
\bq
\small
``\emph{Yes---or no}.''
\eq
\normalsize
\eq

In the economy of no-choice, a contributor, say, a worker or a learner, has only one feasible way of protecting his interests: to silently withhold part of his labour.

The communist block was destroyed by a simple sentiment:

\bq
\small
If they think they pay me let them think I am working.
\normalsize
\eq

Mathematics education in the West is being destroyed by a quiet thought (or even a subconscious impulse):

\bq
\small
If they think they teach me something useful, let them think I am learning.
\eq

On so many occasions I met people who proudly told me:
\bq
\small
I have never been good at mathematics, but I live happily without it.
\normalsize
\eq
They have the right to be proud and confident: they are one-man trade unions who have withheld their learning---and, even they have won nothing, they have not been defeated by the system.

Elizabeth Truss \cite{Truss14} proposes a ``supply-side reform'' of education and skills training as a solution to the hourglass crisis. But supply-side stimuli work best for large scale manufacturers and suppliers. In mathematics education, the key links in the supply chain are children themselves and their families; in the global ``knowledge economy'' too many of them occupy a niche at best similar to that of subsistence farmers in global food production, at worst similar to that of refugees living on food donations. And supply-side economics does not work for subsistence farmers, who, for the escape from the poverty trap, need \emph{demand} for their work and their products, and demand with \emph{payment in advance}---not in 15 or 20 years. Mathematics education has a 15 years long production cycle, which makes supply-side stimuli meaningless.

An additional pressure on mathematics education in the West is created by the division of labour at an international level: in low wage economies of countries like India, learning mathematics still produces economic returns for learners that are sufficiently high in relation to meagre background wages and therefore stimulate ardent learning. As a result, the West is losing the ability to produce competitively educated workers for mathematically intensive industries.

 Should we be surprised that the pyramid of mathematics education is no longer a pyramid and collapses?

\section{The neck of the hourglass}

The mathematical content of the neck can be described in  educationalist terminology used in England as Key Stage 3 (when pupils are aged between 11 and 14) and Key Stage 4 (when pupils are aged between 14 and 16) mathematics:

\bq
\small
Key Stage 3 mathematics teaching [\dots ] marks a transition from the more \emph{informal} approach in primary schools to the formal, \emph{more abstract} mathematics of Key Stage 4 and beyond. \cite[p.~6]{Gardiner14}
\normalsize
\eq
It is informal concrete mathematics and more abstract formal mathematics that make the two bulbs of the hourglass.

Why do we need abstract mathematics? A highly simplified explanation might begin with the fact that money, as it functions in the modern electronic world, is a mathematical abstraction, and this abstraction rules the world.

Of course, this always was the case. However, in 1897 competent handling of money required little beyond arithmetic and the use of tables of compound interest, and clerks at the Post Office were supposed to be mathematically competent for everyday retail finance (see Questions 7 and 8 in the Appendix). Nowadays, the mathematical machinery of finance includes stochastic analysis, among other things. Worse, the mathematics behind the information technology  that supports financial transactions is also very abstract.

Let us slightly scratch the touchscreen of a smartphone or tablet and look at what is hiding behind the ordinary \emph{spreadsheet}.

I prepared the following example for my response to a report from ACME \emph{Mathematical Needs: Mathematics in the workplace
and in Higher Education} \cite{ACME11}\footnote{I used this example in my paper \cite{Borovik12}.}. The report provides the following case study as an important example of use of mathematics.

\bq
\small
\textbf{6.1.4 Case study: Modelling the cost of a sandwich}\\
The food operations controller of a catering company that supplies
sandwiches and lunches both through mobile vans and as special
orders for external customers has developed a spreadsheet that
enables the cost of sandwiches and similar items to be calculated. [\dots ]

\normalsize
\eq

This task would not be too challenging to Post Office clerks of 1897, and would be dealt with by ordinary arithmetic---with the important exception of the ``development of a spreadsheet''. Let us look at it in more detail.

\begin{figure}[h]
\begin{center}
\includegraphics[width=4.5in]{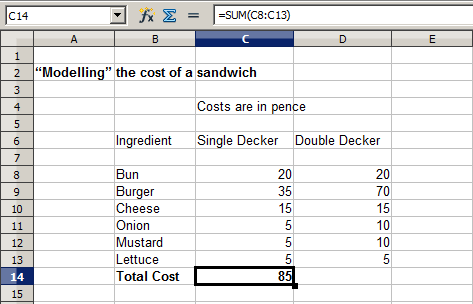}
\end{center}
\caption{}
\label{fig:C14}
\end{figure}

Anyone who ever worked with a spreadsheet of the complexity required for the steps involved in producing sandwiches should know that the key mathematical skill needed is an awareness of the role of brackets in arithmetical expressions and an intuitive feeling for how the brackets are manipulated, something that is sometimes called ``structural arithmetic'' \cite{Gardiner14} or ``pre-algebra''. At a slightly more advanced level working with spreadsheets requires an
understanding of the concept of functional dependency in its algebraic aspects (frequently
ignored in pre-calculus).

To illustrate this point, I prepared a very simple spreadsheet in \textsc{OpenOffice.org Calc} (it uses
essentially the same interface as \textsc{Microsoft Excel}).

\begin{figure}[h]
\begin{center}
\includegraphics[width=4.5in]{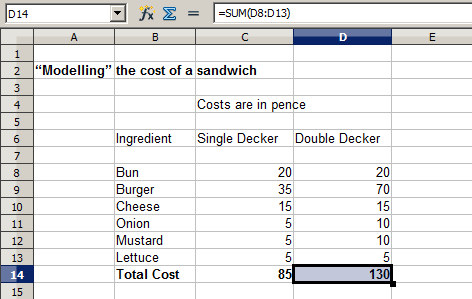}
\end{center}
\caption{}
\label{fig:D14}
\end{figure}

Look at Figure~\ref{fig:C14}: if the content of cell \texttt{C14} is \texttt{SUM(C8:C13)} and you copy cell \texttt{C14} into cell
\texttt{D14} (see Figure~\ref{fig:D14}), the content of cell \texttt{D14} becomes \texttt{SUM(D8:D13)} and thus involves a change of variables. What is
copied is the \emph{structure} of an algebraic expression, not even an algebraic expression itself. And of
course this is not copying the \emph{value} of this expression: notice that the value $85$ becomes
$130$ when moved from cell \texttt{C14} to cell \texttt{D14}!

Intuitive understanding that \texttt{SUM(C8:C13)} is in a sense the same as \texttt{SUM(D8:D13)} can be achieved, for example, by exposing a student to a variety of algebraic problems which convince him/her that a
polynomial of a kind of $x^2 + 2x + 1$ is, from an algebraic point of view, the same as $z^2 + 2z + 1$, and that
in a similar vein, the sum
\begin{center}
\texttt{C8 + C9 + C10 + C11 + C12 + C13}
\end{center}
 is in some sense the same as
\begin{center}
\texttt{D8 + D9 + D10 +
D11 + D12 + D13}.
\end{center}
However the computer programmer (the one who does not merely use spreadsheets, but who writes background code for them), needs an understanding of what it means for two expressions to be ``the same''. Experience suggests rather clearly that the majority of graduates from mathematics departments of British universities, as well as the majority of British school mathematics teachers, do not possess language that allows them to define what  it means for  two expressions in a computer code involving  different symbols (and, frequently, different operations) to be ``actually the same''.

This is a general rule: when a certain previously ``manual'' mathematical procedure is replaced by software, the design and coding of this software requires a much higher level of mathematical skills than is needed for the procedure which has been replaced---but from a much smaller group of workers.

\section{Long division}

For simplistic discussions in the media, the neck of the hourglass  can be summarised in just two words:
\begin{center}
\textbf{long division}.
\end{center}
One of my colleagues who read an early draft of this paper wrote to me:
\bq
\small
``I would not touch long division, as an example, with a ten-foot pole, because it leads to wars.''
\normalsize
\eq

But I am touching it exactly because it leads to wars---to the degree that the words ``long division''  are used  as a symbol for the socio-economic split in English education \cite{Clifton-Crook12}.

Why is long division so divisive? Because it is remarkably useless in the everyday life of 99\% of people. We have to accept that the majority of the population do not need ``practical'' mathematics beyond the use of a calculator, and from the ``practical'' point of view long division can follow slide rules and logarithm tables into the dustbin of history.\footnote{I heard claims that fractions have to be excluded from the school curriculum for the same reason: only a small minority of school students will ever need them in real life.
\bq
``Who of the colleagues present here have lately had to add $\displaystyle{\frac{2}{3}}$ and $\displaystyle{\frac{3}{7}}$?''
\eq
---this question was asked at one of the recent meetings of experts in mathematics education.} But why are long multiplication and long division so critical for squeezing the learners through the hourglass neck?  Because many mathematicians and mathematics educators believe that these ``formal written methods'' should be introduced at a relatively early stage not because of their ``real life relevance'' but with the aim of facilitating children's deep interiorisation of the crucially important class of recursive algorithms which will make the basis of children's later understanding of polynomial algebra---and, at later stages,  ``semi-numerical'' algorithms, in the terminology of the great Donald Knuth \cite{Knuth81}. However there is nothing exceptional about long division: many other algorithms can play in mathematics education the same propaedeutic role, and all of them could be similarly dismissed as not having any ``real life relevance'' because they are needed only by a relatively narrow band of students, those who are expected to continue to learn mathematics up to a more advanced stages and to work in mathematics-intensive industries. In short, ``long division'' is an exemplification of what I later in this paper call ``\emph{deep mathematics education}''.

The left-wing camp in education draw a natural conclusion: long division is hard, its teaching is time- and labour-consuming and therefore expensive, and it will eventually be useful only for a small group of high-flyers---so why bother to teach it?

This is indeed the core question:
\bq
\textbf{Does the nation have to invest human and financial resources into pushing everyone through the hourglass neck? Or should it make a conscious effort to improve the quality of  mathematics teaching, but only for a limited number of students?}
\eq

This is the old conundrum of the British system of education. A recent BBC programme \cite{BBC14} has revealed that Prince Charles in the past lobbied for more (academically selective) grammar schools. The former Education Secretary (Labour) David Blunkett told about his exchanges with Prince Charles:

\bq
\small
I would explain that our policy was not to expand grammar schools, and he didn't like that.

He was very keen that we should go back to a different era where youngsters had what he would have seen as \textbf{the opportunity to escape from their background, whereas I wanted to change their background}. [The emphasis is mine---AB.]
\normalsize
\eq

This is a brilliant formulation of the dilemma, and it is especially good in the case of mathematics education because the hourglass shape of economic demand for different levels of mathematics education puts the emotive word ``escape'' on a solid economic foundation: it is the escape through the hourglass neck.

While I would be delighted, and relieved, to be convinced by arguments to the contrary, at this point I can see the solutions offered by the Left and the Right of British education politics as deficient in ways that mirror each other:
\bi
\item The Left appear to claim that it is possible to have quality mathematics education for everyone. While their position is sincerely held, still, as I see it, it leads to inconsistencies  which can be avoided only by lowering the benchmark of ``quality'' and ignoring the simple economic fact that what they call ``quality education'' is neither needed by, nor required from, learners in their life, present and future, outside school.

\item The Right appear to claim that administrative enforcement of standards will automatically raise the quality of education for everyone. It is also a sincerely held position, but, as  I see it, it leads to inconsistencies  which can be avoided only by preparing escape routes for their potential voters' children in the form of ``free schools''.
\ei
My previous analysis has not made any distinction between ``state'' and ``private'' schools; this reflects my position---I do not believe that mainstream private schools, or ``free schools'' (even it they are privatised in the future) make any difference in the systemic crisis of mathematics education.

\section{Back to \emph{Z\"{u}nfte}?}

In relation to mathematics, social factors and, consequently, social division of labour attain increasing importance for a simple reason:  who but families are prepared to invest 15 years into something so increasingly specialised as mathematics education?

What instructional system was in place before the division-of-labor sweatshops glorified by Adam Smith? The \emph{Zunft system}. In German, \emph{Zunft} is a historic term for a guild of master craftsmen (as opposed to trade guilds). The high level of specialisation of \emph{Z\"{u}nfte} could be sustained only by hereditary membership and training of craftsmen, from an early age, often in a family setting. It is hard not to notice a certain historical irony\dots

The changing patterns of division of labour affect mathematics education in every country in the world. But reactions of the government, of the education community, of parents from different social strata depend on the political and economic environment of every specific country. So far I analysed consequences for education policy in England; when looking overseas and beyond the anglophone world, one of more interesting  trends is mathematics homeschooling and ``math circles'' movements in two countries so different as the USA and Russia. In both countries mathematically educated sections of middle class are losing confidence in their governments' education policies and in the competence of the mathematics education establishment, and are choosing to pass on their own expertise through homeschooling as a modern \emph{Zunft}.

Some of the economic forces affecting education are brutally simple, and the principal barrier facing potential homeschoolers is purely financial. Mainstream education fulfils an important function of a storage room for children, releasing parents for salaried jobs; if parents were to spend more time with children, rates of pay would have to be higher. A family cannot homeschool their children without sufficient disposable income, part of which can be re-directed and converted into ``quality time'' with children.

Statistics of mathematics homeschooling are elusive, but what is obvious is the highest quality of intellectual effort invested in the movement by its leading activists---just have a look at books \cite{Burago12,Droujkova14,Zvonkin11}. At the didactical level, many  inventions of mathematics homeschoolers are wonderful but intrinsically unscalable and cannot be transplanted into the existing system of mass education. I would say that their approach is not a remedy for the maladies of mainstream education; on the contrary, the very existence of mathematics homeschoolers is a symptom of, and a basis for a not very  optimistic prognosis for, the state of mass mathematics education.

Still, in my opinion, no-one in the West has captured the essence of \emph{deep mathematics education} better then they have.

\section{\emph{Z\"{u}nfte} and ``deep mathematics education''}

At the didactic level,  bypassing the hourglass neck of economic demand for mathematics means development of \emph{deep mathematics education}. I would define it as
\bq\small
Mathematics education in which every stage, starting from pre-school, is designed to fit the individual cognitive profile of the child and to serve as propaedeutics of his/her even deeper study of mathematics at later stages of education---including transition to higher level of abstraction and changes of conceptual frameworks.
\normalsize
\eq
To meet these aims, ``deep'' mathematics education should unavoidably be joined-up and cohesive.\footnote{The Moscow Center for Continuous   Mathematics Education, \url{http://www.mccme.ru/index-e1.html}, emphasises this aspect by putting the word ``continuous'' into its name; it focuses on bridging the gap between school and university level mathematics, while homeschoolers tend to start at the pre-school stage.}

To give a small example in addition to the already discussed long division, I use another stumbling block of the English National Curriculum: times tables. The following is a statutory requirement:

\bq\small
By the end of year 4, pupils should have memorised their multiplication tables up to and including the 12 multiplication table and show precision and fluency in their work. \cite{DfE13}
\normalsize \eq
This requirement is much criticised for being archaic (indeed, why $12$?), cruel and unnecessary. But to pass through the neck of the hourglass, children should know by heart times tables up to 9 by 9; even more, it is very desirable that they know by heart square numbers up to $20^2 = 400$, because understanding and ``intuitive feel'' of behaviour of quadratic functions is critically important for learning algebra and elementary calculus.

The concept of ``deep mathematics education'' is not my invention.
I borrowed the words  from Maria Droujkova, one of the leaders of mathematics homeschooling. Her understanding of this term is, first of all, deeply human and holistic.

In her own words\footnote{Private communication.},

\bq\small
 The math we do is defined by freedom and making. We value mastery---with the understanding that different people will choose to reach different levels of it. The stances of freedom and making are in the company's motto:
\bq \emph{Make math your own, to make your own math.}
\eq

When I use the word ``deep'' as applied to mathematics education, I approach it from that natural math angle. It means deep agency and autonomy of all participants, leading to deep personal and communal meaning and significance; as a corollary, deep individualization of every person's path; and deep psychological and technological tools to support these paths.
\normalsize\eq

Droujkova uses, as an example, iterative algorithms, and her approach to this concept is highly relevant for the discussion of the propaedeutic role of ``long division'':
\bq\small

From the time they are toddlers, children play with recursion and iteration, in the contexts where they can define their own iterating actions. For example, children design input-output ``function machines'' and connect the output back to the input. Or experiment with iterative splitting, folding, doubling, cutting with paper, modeling clay, or virtual models. Or come up with substitution and tree fractals, building several levels of the structure by iterating an element and a transformation. Grown-ups help children notice the commonalities between these different activities, help children develop the vocabulary of recursive and iterative algorithms, and support noticing, tweaking, remixing, and formulating of particular properties and patterns. As children mature, their focus shifts from making and remixing individual algorithms to purposeful creation and meta-analysis of patterns. For example, at that level children can compare and contrast recursion and iteration, or analyze information-processing aspects of why people find recursive structures beautiful, or research optimization of a class of recursive algorithms.
\normalsize\eq

Maria Droujkova describes a rich and exciting learning activity. But it would be impossible without full and informed support from children's families. To bring this education programme to life, you need a community of like-minded and well-educated parents. It could form around their children's school (and would almost inevitably attempt to control the school), or around a ``mathematical circle'', informal and invisible to the educational establishment and therefore free from administrative interference; or, what is much more likely in our information technology age, it could grow as an Internet-based network of local circles connected by efficient communications tools---and perhaps helped by parents' networking in their professional spheres. These ``\emph{communities of practice}'', as Droujkova calls them using a term coined by Wenger \cite{Wenger00}, are \emph{Z\"{u}nfte} at the new turn of history's spiral. I see nothing that makes them unfeasible.

I wish  mathematics homeschoolers the best of luck. But their work is not a recipe  for mainstream education.

\section{``Deep mathematics education'': Education vs.\ training}

\medskip
\small
\begin{flushright}
\emph{
Who knows the difference between education and training?\\ For those of you with daughters, would you rather have\\ them take sex education or sex training? Need I say more?}\\ Dennis Rubin\\[2ex]
\end{flushright}
\normalsize

The witticism above
makes it clear what is expected from education as opposed to training: the former should give a student ability to make informed and responsible \emph{decisions}.

The same message is contained in the apocryphal saying traditionally attributed to a  President of Harvard University who  allegedly said,  in response to a question on what was so special about Harvard to justify the extortionate fees,

\bq
\emph{``We teach criteria.''}
\eq

Let us think a bit: who needs criteria? Apparently, people who, in their lives,  have to make choices and decisions. But millions of people around us are not given the luxury of choice.

This is the old class divide that tears many education systems apart: education is for people who expect to give orders; training is for ones who take orders.
Mathematics, as it is taught in many schools and universities, is frequently reduced to training in a specific transferable skill: the ability to carry out meaningless repetitive tasks. Unfortunately, many of the students who I meet in my professional life have been, in my assessment, trained, not educated: they have been taught to the test, and at the level of rudimentary procedural skills which can be described as a kind of  painting-by-numbers.

\begin{figure}[h]
\begin{center}
\includegraphics[width=4.5in]{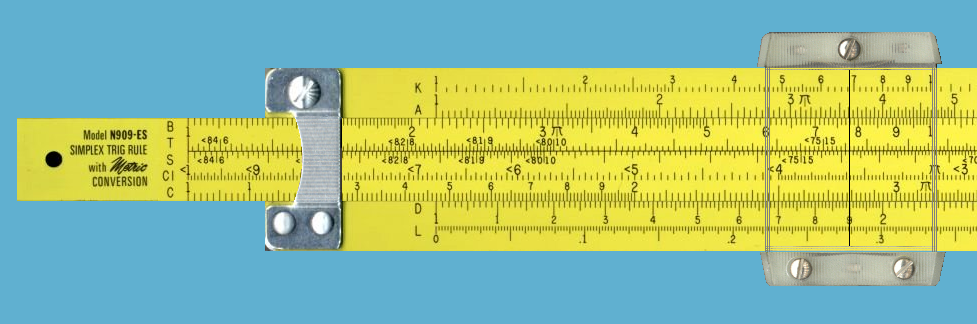}
\end{center}
\caption{
Simulated Pickett N909-ES Slide Rule. It is fully functional (but needs a sufficiently wide computer screen)! Source: {%
http://www.antiquark.com/sliderule/sim/n909es/virtual-n909-es.html}.}
\label{fig:sliderule}
\end{figure}

This divide between education and training remains a forbidden theme in mathematics education discourse in England. But a better understanding of what makes education different from training would help, for example, in the assessment of possibilities offered by new computer-assisted and computer-based approaches to mathematics learning and teaching. I would not be surprised if computerisation of mathematics training could be achieved easily and on the cheap---but I also think that any attempt to do that is likely to be self-defeating.  Indeed I believe in a basic guiding principle: if a certain mathematical skill can be taught by a computer, this is the best proof that this skill is economically redundant---it could be best done by computers without human participation, or outsourced to a country with cheaper labour. (For readers who remember slide rules, this is like using computers for teaching and learning the technique of slide rule calculations. By the way, you can find on the Internet fully functional virtual slide rules, with moving bits that can be dragged by a mouse, see Figure~\ref{fig:sliderule}.)

\begin{figure}[h]
\begin{center}
\includegraphics[width=4.5in]{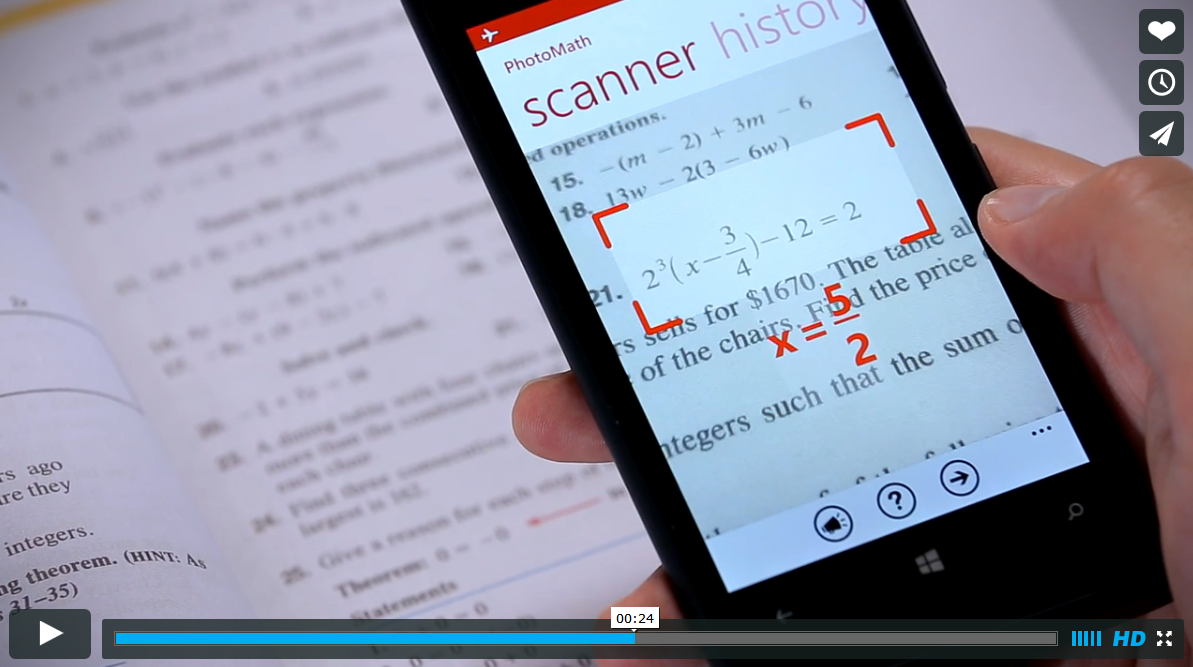}
\end{center}
\caption{A screen shot from an advert for \textsc{PhotoMath}: point your smartphone at a problem in the textbook, and the answer is instantly produced. Source: {http://vimeo.com/109405701}.}
\label{fig:PhotoMath}
\end{figure}

\begin{figure}[h]
\begin{center}
\includegraphics[width=4.5in]{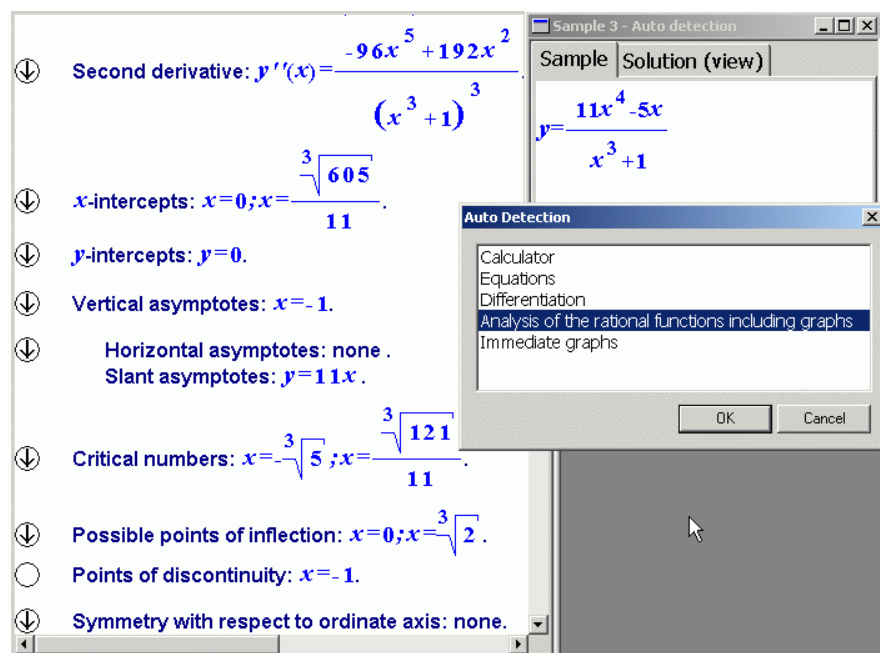}
\end{center}
\caption{A screen shot from \textsc{Universal Math Solver}: a few intermediate steps of the analysis of behaviour of the function
$\displaystyle{y = \frac{11x^4 - 5x}{x^3+1}}$.
``Arrow Down'' icons on the left margin unroll a more complete write-up for particular steps in calculations. Figure~\ref{fig:MathsolverB} shows the graph of the function.  Source: {http://www.umsolver.com/}.}
\label{fig:Mathsolver}
\end{figure}

\begin{figure}[h]
\begin{center}
\includegraphics[width=4.5in]{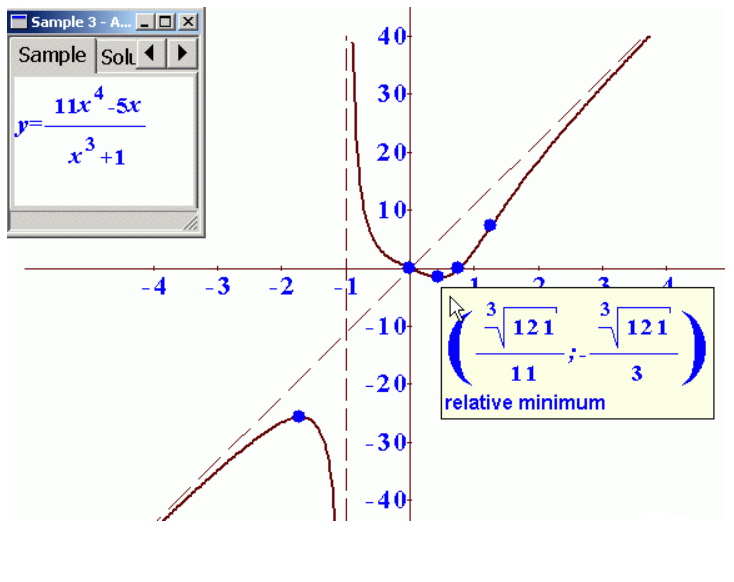}
\end{center}
\caption{A screen shot from \textsc{Universal Math Solver}: a graph of the function $\displaystyle{y = \frac{11x^4 - 5x}{x^3+1}}$
which highlights information obtained at the previous steps of analysis, see Figure~\ref{fig:Mathsolver}.  Source: {http://www.umsolver.com/}.}
\label{fig:MathsolverB}
\end{figure}

Unfortunately, almost the entire school and a significant part of undergraduate mathematics, as it is currently taught in England, is likely to follow the slide rules into the dustbins of history.
Figure~\ref{fig:PhotoMath} shows an advert for a smartphone app \textsc{PhotoMath}, it has gone viral and enjoys an enthusiastic welcome on the Internet. Mathematical capabilities of \textsc{PhotoMath}, judging by the product website\footnote{\url{%
http://www.windowsphone.com/en-us/store/app/photomath/1f25d5bd-9e38-43f2-a507} \url{-a8bccc36f2e6}.} are still relatively modest. However, if the scanning and optical character recognition modules of \textsc{PhotoMath} are combined with the full version of Yuri Matiasevich's \textsc{Universal Math Solver}, it will solve at once any mathematical equation or inequality,  or evaluate any integral, or check convergence of any series appearing in the British school and undergraduate mathematics. Moreover, it will produce, at a level of detail that can be chosen by the user, a complete write-up of a solution, with all its cases, sub-cases, and necessary explanations. Figures~\ref{fig:Mathsolver} and \ref{fig:MathsolverB} show that, unlike industrial strength software packages \textsc{Maple} and \textsc{Mathematica},  \textsc{Universal Math Solver} faithfully follows the classical ``manual'' procedures of mathematics textbooks.

This presents a historically  unprecedented challenge to the teaching profession: how are we supposed to teach mathematics to students who, from age of five, have on their smartphones, or on smartglasses, or other kinds of wearable smart devices, apps that instantly answer every question and solve every problem from school and university textbooks?

In short, smart phones can do exams better than humans, and the system of ``procedural'' mathematics training underpinned by standardised written examinations is dead. Perhaps, we have to wait a few years for a coroner's report, but we can no longer pretend that nothing has happened.

By contrast, ``deep mathematics education''   treats mathematics as a discipline and art of those aspects of formal reasoning \emph{which cannot be entrusted to a computer}. This is, in  essence, what mathematics homeschoolers are trying to develop.

I am a bit more cautious about the feasibility of setting-up and developing a system of ``deep mathematics education''  at a national level. It is likely to be expensive and raises a number of uncomfortable political questions. To give just one example of a relatively benign kind: in such a system, it could be desirable to have oral examinations in place of written ones. The reader familiar with the British university system, for example, can easily imagine all the political complications that would follow.

\section{``Deep mathematics education'': Phase transitions and metamorphoses}

\medskip
\small
\begin{flushright}
\emph{We are caterpillars of angels.}\\
Vladimir Nabokov
\end{flushright}
\normalsize

I am old enough to have been taught, in my teenage years, to write computer code in physical addresses, that is, sequences of zeroes and ones, each sequence referring to a particular memory cell in the computer. My colleague, an IT expert, told me recently that he and people who work for him passed in their lives through 6 (six!) changes of paradigms of computer programming. In many walks of life, to have a happy and satisfying professional career, one has to be future-proof  by being able to re-learn the craft, to change his/her way of thinking.

How can this skill of changing one's way of thinking be acquired and nurtured?

At school level---mostly by learning mathematics. Regular and unavoidable changes of mathematical language reflect changes of mathematical thinking. This makes mathematics different from the majority of other disciplines.

The crystallisation of a mathematical concept (say, of a fraction), in a child's mind could be like a phase transition in a crystal growing in a rich, saturated---and undisturbed---solution of salt. An ``aha!'' moment is a sudden jump to another level of abstraction. Such changes in one's mode of thinking are like a metamorphosis of a caterpillar into a butterfly.

As a rule, the difficulties of learning mathematics are difficulties of adjusting to change. Pupils who have gained experience of overcoming these difficulties are more likely to grow up future-proof. I lived through sufficiently many changes in technology to become convinced that mathematically educated people are stem cells of a technologically advanced society, they are re-educable, they have a capacity for metamorphosis.

As an example of a sequence of paradigm changes in the process of learning, consider one of the possible paths in learning algebra.
I picked this path because it involves  three ``advanced'' concepts which, in the opinion of some educationalists, can be removed from mainstream school mathematics education as something that has no practical value: fractions and long division (which featured earlier in this paper), and factorisation of polynomials.

\begin{figure}
\begin{center}
\includegraphics[width=3in]{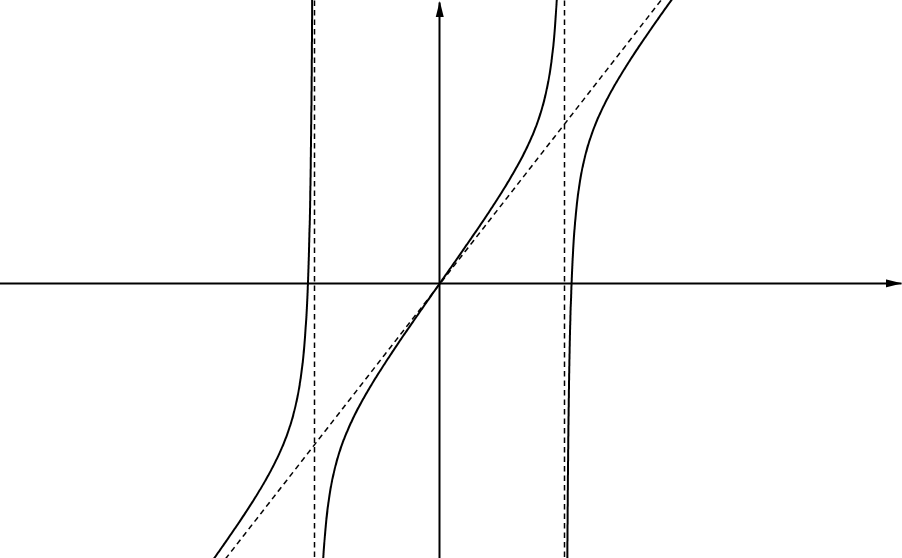}
\end{center}
\caption{This problem:
 \emph{``Find a rational function which has a graph with vertical and oblique asymptotes as shown on this drawing''}
is a long way from the primary school fractions and ratios, but it is about ratios---this time of variable quantities. It is even more useless in the ``everyday life'' than fractions; its value lies in providing an example of a link between algebra, geometry and topology as well as giving a tangible example of ``asymptotic behaviour'', a concept of crucial importance for many applications of mathematics.  Perhaps someone who has not mastered fractions at primary school still has a chance to reach, in his/her later years, the level of understanding of elementary algebra and pre-calculus necessary for solving this problem, but this is likely to be the exception rather than the rule. The problem and drawing by Julian Gilbey, reproduced with his kind permission.}
\label{fig:obliqueasymptote}
\end{figure}

The path, one of many in mathematics learning, goes from pre-school to undergraduate courses:
\begin{enumerate}
\item   Naive arithmetic of natural numbers;

\item fractions and negative numbers;

\item place value, formal written algorithms (the ``long multiplication'' is the most important of them), ``structural arithmetic'' (that is, ability to simplify arithmetic calculations such as  $17 \times 5 + 3 \times 5$);

\item algebraic notation;

\item polynomials; roots and  factorisation of polynomials as a way to see that polynomials \emph{have their own life} in a new mathematical world, much wider and richer than arithmetic---in particular, this means that ``long multiplication'' and ``long division'' are revisited in symbolic form;

\item interpretation of polynomials as functions;  coordinates and graphs;

\item rational functions (ratios of polynomials) in two facets: as fractions revisited in symbolic form, and as  functions;

\item and something that is not usually mentioned in school mathematics: understanding that the behaviour of a rational function $f(x)/g(x)$ \emph{as a function} is dictated by its zeroes and poles (singularities), that is, by roots of the numerator $f(x)$ and denominator $g(x)$, thus revisiting factorisation at a new level---see Figure~\ref{fig:obliqueasymptote} for an example;

\item and, finally, something that is not always mentioned in undergraduate courses: the convergence radii of the power series
\[
\frac{1}{1+x^2} = 1 -x^2 + x^4 - x^6 + \cdots
\]
and
\[
\arctan x = x - \frac{x^3}{3} + \frac{x^5}{5} - \frac{x^7}{7} + \cdots
\]
equal $1$ because, in the complex domain, the first of the two functions, the rational (and hence analytic) function
\[
f(z) = \frac{1}{1+z^2}
\]
has poles  $z = i$ and $z = -i$, both at distance $1$ from $0$, and because the second function is an integral of the first one
\[
\arctan z = \int \frac{dz}{1+z^2}.
\]
\end{enumerate}
Even ignoring the stages 8 and 9, we have six deep and difficult changes of the mathematical language used and of the way of thinking about mathematical objects. Each of these six steps is challenging for the learner. But  they constitute a good preparation for facing and overcoming  future changes in professional work.

I have used the classical school algebra course and a bit of calculus as an example. I accept that mathematics can be taught differently.  I myself can offer some modifications---for example, why not introduce children, somewhere after level 1, to a toy object-oriented programming language of the kind of \textsc{ScratchJr}\footnote{\textsc{ScratchJr} allows the learner  to build iterative algorithms---see a discussion of their pedagogical value in  Droujkova's quote above---by moving and snapping together \textsc{Lego}-style blocks on a touchscreen.}, and, after level 6, to some appropriately simplified version of a  \textsc{Haskell}-kind \cite{Hutton07} language of functional programming? But, I wish to re-iterate, I refrain from any recommendations, especially if they require a mass scale re-education of the army of teachers.

However, every approach to learning mathematics, if it leads to a certain level of mastering mathematics, will inevitably involve several changes of the underlying conceptual framework and the language of mathematical expression, at every stage increasing the level of abstraction and the compression of information. What eventually matters is the degree of compression (and the latter more or less correlates with the number of phases of development  through which a student passed). Many undergraduate mathematics students come to university with a depleted ability to compress their mathematical language further, and this is happening because their  previous ``phase transitions'' were badly handled by their teachers.

\bq
\textbf{The potential for further intellectual metamorphoses is the most precious gift of ``deep mathematics education''.}
\eq

\section{Conclusion}

\medskip

\small
\begin{flushright}
\emph{I came here knowing we have some sickness in our system\\ of education; what I have learned is that we have a \emph{\textbf{cancer!}}}\\
Richard Feynman, \emph{Surely You're Joking, Mr. Feynman!}\\[2ex]
\end{flushright}
\normalsize

In this paper, I have attempted to describe how deepening specialisation and division of labour in the economy affects the  mathematics education system, changes its shape, undermines its stability, leads to a social split in mathematics education, and (at least in England) provokes political infighting.

I wish to reiterate that I am not taking sides in these fights. I do not wish to lay blame on anyone, or criticise anyone's views. My paper is a call for a sober, calm, and apolitical discussion of the socio-economic roots of the current crisis in mathematics education.

Mathematics at the level needed for serious work, say,  in electronics and information technology,  requires at least 15 years of
systematic stage-by-stage learning, where steps cannot be arbitrarily swapped or skipped. After all, it's about growing neuron connections in the brain, it is a slow process. Also, it is an age-specific process, like learning languages.

Democratic nations, if they are sufficiently wealthy,  have three options:
\bi

\item[\hspace{-1em}(A)] Avoid limiting children's future choices of profession, teach rich mathematics to every child---and invest serious money into thorough
professional education and development of teachers.

\item[\hspace{-1em}(B)] Teach proper mathematics, and from an early age, but only to a selected minority of children. This is a much
cheaper option, and it still meets the requirements of industry, defence and security sectors, etc.

\item[\hspace{-1em}(C)] Do not teach proper mathematics at all and depend on other countries for the supply of technology  and  military
protection.
\ei

Which of these options are realistic in a particular country at a given time, and what the choice should be, is for others to decide.

I am only calling a spade a spade.

\section*{Acknowledgements}

I first used the ``pyramid'' diagram in my talk at the Mathematics Colloquium at the Middle East Technical University, Turkey, in April 2007, and I thank Ay\c{s}e Berkman for her kind invitation to give that talk. The paper was developed into its present form thanks to my involvement with CMEP, the Cambridge Mathematics Education Project. I thank my CMEP colleagues for many useful discussions---but neither they nor CMEP are responsible for my views expressed here.

I am deeply grateful to Julia Brodsky, Dmitri Droujkov and Maria Droujkova for generously sharing their ideas.

I thank Frank Wagner for finding for me the word \emph{Zunft}; I knew its Russian analogues, but could not find an appropriate English term.

Michael Barany, Gregory Cherlin, David Edwards, Rosemary Emmanuel, Jean-Michel Kantor, Alexander Kheyfits, Roman Kossak, David Pierce, Seb Schmoller, Victor Sirotin, and several mathematician colleagues who preferred to stay anonymous provided useful comments on my text and/or helped to improve its style and grammar.  I am grateful to them all---but they do not bear any responsibility for what is said in this paper.

Julian Gilbey kindly allowed to reproduce his problem and Figure~\ref{fig:obliqueasymptote}.

And I thank the anonymous referee for the most helpful advice.

\newpage
\small

\begin{center}
{\Large\textsc{Appendix}}\\
\textsc{Post Office Entrance Examination\\
Women And Girl Clerks}\\
October 1897
\end{center}

\fff{1} Simplify
\[
\frac{1/2+1/3+1/4+1/5}{1/2+1/3-1/4-1/5}
+ \frac{1/4+1/5+1/6+1/7}{1/4+1/5-1/6-1/7}
- \frac{1024}{1357} .
\]

\fff{2} If $725$ tons $11$ cwts. $3$ qrs. $17$ lbs. of potatoes cost \pounds $3386$, $2$s. $2\frac{1}{2}$d. how much will 25 tons 11 cwts.
3 qrs. 17 lbs. costs (sic)?

\fff{3} Extract the square root of $331930385956$.

\fff{4} A purse contains 43 foreign coins, the value of each of which either exceeds or falls short of one crown
by the same integral number of pence. If the whole contents of the purse are worth \pounds $10$, $14$s. $7$d., find the
value and number of each kind of coin. Show that there are two solutions.

\fff{5} Explain on what principle you determine the order of the operations in
\[
\frac{1}{2} + \frac{3}{4} \div \frac{5}{6} - \frac{7}{8} \times \frac{9}{10} ,
\]
and
express the value as a decimal fraction. Insert the brackets necessary to make the expression mean :-

\begin{quote}
Add $\frac{3}{4}$ to $\frac{1}{2}$, divide the sum by $\frac{5}{6}$, from the quotient subtract $\frac{7}{8}$, and multiply this difference by $\frac{9}{10}$.
\end{quote}

\fff{6} Show that the more figures 2 there are in the fraction $0.222\dots 2$, the nearer its value is to $\frac{2}{9}$. Find the difference in value when there are ten $2$'s.

\fff{7} I purchased \pounds $600$ worth of Indian 3 per cent. stock at 120. How much Canadian 5 per cent. stock
at 150 must I purchase in order to gain an average interest of 3 per cent. on the two investments (sic!)?

\fff{8} If five men complete all but $156$ yards of a certain railway embankment, and seven men could
complete all but $50$ yards of the same embankment at the same time, find the length of the embankment.

\fff{9} Find, to the nearest day, how long \pounds $390$, $17$s. $1$d. will take to amount to \pounds $405$, $14$s. $3$d. at $3\frac{1}{4}$ per
cent. per annum ($365$ days) simple interest.

\fff{10} A certain Irish village which once contained $230$ inhabitants, has since lost by emigration three-fourths of its agricultural population and also five other inhabitants. If the agricultural population is now as numerous as the rest, find how the population was originally divided.

\


\begin{thebibliography}{99}

\bibitem{abel-deitz} J. R. Abel and R. Deitz, Job polarization in the United States: a widening gap and shrinking middle,
 Liberty Street Economics, 21 November 2011, \url{http://libertystreeteconomics.newyorkfed.org/2011/11/job-polarization-in-the-united-states-a-widening-gap-and-shrinking-middle.html}. Accessed 11 July 2012.

\bibitem{Acemoglu10}
D. Acemoglu and D. Autor,
Skills, tasks and technologies: implications for employment and
earnings, 10 June 2010. \url{http://economics.mit.edu/files/5571}, accessed 11 July 2012.

\bibitem{ACME11} ACME, Mathematical needs. Mathematics in the workplace
and in higher education. Advisory Committy on Mathematics Education (ACME), June 2011.
ISBN 978-0-85403-895-4.



\bibitem{Anderson09}
P. Anderson, Intermediate occupations and the conceptual and empirical limitations of the hourglass economy thesis,
Work, employment and society 23 no.\ 1(2009), 169--180. DOI: 10.1177/0950017008099785.



\bibitem{Autor10}
D. Autor, The polarization of job opportunities in the U.S. labor market: implications for employment and earnings. Center for American Progress and The Hamilton Project, May 2010. \url{http://economics.mit.edu/files/5554}, accessed 11 July 2012.


\bibitem{BBC14} BBC,  Prince Charles `tried to influence government decisions'. 29 June 2014. \url{http://www.bbc.com/news/uk-politics-28066081}.

\bibitem{Borovik12} A. Borovik, L'exemple de la modelisation, Commentaire
138 (2012) 500--502.



\bibitem{Braverman74} H. Braverman, Labor and monopoly capital: the degradation of work in the twentieth century. 1974. Republished: Monthly Review Press, 1998. ISBN: 9780853459408

\bibitem{Brynjolfsson14}    E. Brynjolfsson and A. McAfee, The second machine age: work, progress, and prosperity in a time of brilliant technologies. W. W. Norton \&\ Company, 2014. ISBN-10: 0393239357,
ISBN-13: 978-0393239355.

\bibitem{Burago12} A. Burago,
Mathematical circle diaries, Year 1: complete curriculum for grades 5 to 7. American Mathematical Society, 2012.
ISBN-10: 0821887459.
ISBN-13: 978-0821887455.

\bibitem{Clifton-Crook12} J. Clifton and W. Cook, A long division: Closing the attainment gap in England's
secondary schools. Institute for Public Policy Research, 2012. \url{http://www.ippr.org/assets/media/images/media/files/publication/2012/09/long division FINAL version_9585.pdf}.


\bibitem{Cowen13} T. Cowen, Average is over: oowering America beyond the age of the great stagnation. Penguin, 2013. ISBN-13:9780525953739. ISBN-10:0525953736.


\bibitem{Edwards14} D. Edwards, The math myth. The De Morgan Gazette 5 no. 3 (2014), 19--21. ISSN 2053--1451.
\url{http://education.lms.ac.uk/wp-content/uploads/2014/07/DMG_5_no_3_2014.pdf}.


\bibitem{Elliott11}
P. Elliott, a blog post on 9 July 2011 on the Times Education Supplement website, \url{http://www.timeshighereducation.co.uk/story.asp?sectioncode=26&storycode=416765&c=1}.

\bibitem{Gainsburg05} J. Gainsburg, School mathematics in work and life: what we know
and how we can learn more. Technology in Society 27 (2005) 1--22.

\bibitem{Gardiner14} A. D. Gardiner, Teaching mathematics 	
at secondary level.
The De Morgan Gazette 6 no.~1 (2014), 1--215. ISSN 2053--1451.


\bibitem{Goos09}
M. Goos, A. Manning, and A. Salomons,
Job polarization in Europe.
American Economic Review: Papers \&\ Proceedings  99 no. 2 (2009), 58--63. \url{http://www.econ.kuleuven.ac.be/public/n06022/aerpp09.pdf}.

\bibitem{Greenspan13}
A. Greenspan, The map and the territory: risk, human nature and the future of forecasting. Allen Lane, 2013. ISBN-13: 978 0 241 00359 6.



\bibitem{Hutton07}
G. Hutton,  Programming in Haskell. Cambridge: Cambridge University Press, 2007.

\bibitem{Knuth81} D. E. Knuth, The art of computer programming. {V}ol. 2,
Seminumerical algorithms.
Addison-Wesley, 1981.      ISBN 0-201-03822-6.

\bibitem{Levy04}
F. Levy and R. J. Murnane, The new division of labor. How computers are creating the next job market. Princeton University Press, 2004. ISBN-10: 0-691-12402-7. ISBN-13: 978-0-691-12402-5.

\bibitem{DfE13} Department for Education. Statutory guidance. National curriculum in England: mathematics programmes of study.  2013.
\url{https://www.gov.uk/government/publications/national-curriculum-in-england-mathematics-programmes-of-study/national-curriculum-in-england-mathematics-programmes-of-study#year-6-programme-of-study}.

\bibitem{Droujkova14}
Y. McManaman, M. Droujkova, and E. Salazar,
 Moebius noodles: adventurous math for the playground crowd.  Delta Stream Media, 2014. [Kindle Edition, sold by Amazon.]



\bibitem{Truss11}
E. Truss, Academic rigour and social mobility: how low income students are being kept out of top jobs. Centre:Forum, 2011. \url{http://www.centreforum.org/assets/pubs/academic-rigour-and-social-mobility.pdf}. Accessed 10 July 2011.


\bibitem{Truss14}
E. Truss, The global education race. Speech at
    The Oxford Conference in Education, St John's College, Oxford, 3 January 2014.
\url{https://www.gov.uk/government/speeches/elizabeth-truss-the-global-education-race}, accessed 16 June 2014.



\bibitem{Wenger00} E. Wenger, Communities of practice: learning, meaning, and identity.
Cambridge University Press, 2000.
ISBN-10: 0521663636. ISBN-13: 978-0521663632.


\bibitem{Zvonkin11} A. Zvonkin, Math from three to seven: the story of a mathematical circle for preschoolers. American Mathematical Society, 2011. ISBN-10: 082186873X.
ISBN-13: 978-0821868737.

\end{thebibliography}
\end{document}